\documentclass[a4paper]{article}
\usepackage[left=3.0cm,right=3.0cm,top=2.5cm,bottom=2.5cm]{geometry}
\usepackage{xcolor}
\usepackage{graphicx}
\usepackage{cite}
\usepackage{amsmath,amssymb,amscd,amsxtra,amsfonts}

\newcommand\keywordsname{Key words}
\newcommand\AMSname{AMS subject classifications}

\newenvironment{@abssec}[1]
{\if@twocolumn
\section*{#1}%
\else
\vspace{.05in}\footnotesize
\parindent .2in
{\upshape\bfseries #1. }\ignorespaces
\fi}

{\if@twocolumn\else\par\vspace{.1in}\fi}

\newenvironment{keywords}{\begin{@abssec}{\keywordsname}}{\end{@abssec}}

\providecommand{\Div}{\operatorname{div}}          
\providecommand{\curl}{\operatorname{{\bf curl}}}  

\providecommand*{\Dist}[2]{\operatorname{dist}({#1};{#2})}   
\providecommand*{\Dist}[2]{\Dist{#1}{#2}}

\newcommand{\Ve}{{\mathbf{e}}}

\newcommand{\Vr}{{\mathbf{r}}}

\newcommand{\Vx}{{\mathbf{x}}}

\newcommand{\Bf}{{\boldsymbol{f}}}
\newcommand{\Bg}{{\boldsymbol{g}}}

\newcommand{\Bn}{{\boldsymbol{n}}}

\newcommand{\Bu}{{\boldsymbol{u}}}
\newcommand{\Bv}{{\boldsymbol{v}}}
\newcommand{\Bw}{{\boldsymbol{w}}}
\newcommand{\Bx}{{\boldsymbol{x}}}

\newcommand{\VI}{{\mathbf{I}}}

\newcommand{\BB}{{\boldsymbol{B}}}
\newcommand{\BC}{{\boldsymbol{C}}}

\newcommand{\BH}{{\boldsymbol{H}}}

\newcommand{\BL}{{\boldsymbol{L}}}

\newcommand{\BP}{{\boldsymbol{P}}}

\newcommand{\BV}{{\boldsymbol{V}}}

\newcommand{\varphibf}{\boldsymbol{\varphi}}

\newcommand{\psibf}{\boldsymbol{\psi}}

\newcommand{\Ca}{\mathcal{A}}
\newcommand{\Cb}{\mathcal{B}}

\newcommand{\Cf}{\mathcal{F}}

\newcommand{\Co}{\mathcal{O}}

\newcommand{\Ct}{\mathcal{T}}

\newcommand{\bbA}{\mathbb{A}}
\newcommand{\bbB}{\mathbb{B}}
\newcommand{\bbC}{\mathbb{C}}

\newcommand{\bbF}{\mathbb{F}}
\newcommand{\bbG}{\mathbb{G}}

\newcommand{\bbJ}{\mathbb{J}}

\newcommand{\bbL}{\mathbb{L}}
\newcommand{\bbM}{\mathbb{M}}

\newcommand{\bbP}{\mathbb{P}}
\newcommand{\bbQ}{\mathbb{Q}}
\newcommand{\bbR}{\mathbb{R}}
\newcommand{\bbS}{\mathbb{S}}

\newcommand{\bbX}{\mathbb{X}}

\newcommand*{\N}[1]{\left\|{#1}\right\|}     
\newcommand*{\SN}[1]{\left|{#1}\right|}      

\newcommand*{\Lp}[2][\defaultdomain]{L^{#2}({#1})}
\newcommand*{\Lpv}[2][\defaultdomain]{\BL^{#2}({#1})}
\newcommand*{\NLp}[3][\defaultdomain]{\N{#2}_{\Lp[#1]{#3}}}

\newcommand*{\Ltwo}[1][\defaultdomain]{\Lp[#1]{2}}
\newcommand*{\Ltwov}[1][\defaultdomain]{\Lpv[#1]{2}}
\newcommand*{\NLtwo}[2][\defaultdomain]{\NLp[#1]{#2}{2}}

\newcommand*{\Hm}[2][\defaultdomain]{H^{#2}({#1})}
\newcommand*{\Hmv}[2][\defaultdomain]{\BH^{#2}({#1})}
\newcommand*{\bHm}[3][\defaultdomain]{H_{#3}^{#2}({#1})}

\newcommand*{\bHmv}[3][\defaultdomain]{\BH_{#3}^{#2}({#1})}

\newcommand*{\Hone}[1][\defaultdomain]{\Hm[#1]{1}}
\newcommand*{\Honev}[1][\defaultdomain]{\Hmv[#1]{1}}

\newcommand*{\zbHone}[1][\defaultdomain]{\bHm[#1]{1}{0}}
\newcommand*{\zbHonev}[1][\defaultdomain]{\bHmv[#1]{1}{0}}
\newcommand*{\NHone}[2][\defaultdomain]{{\N{#2}}_{\Hone[{#1}]}}
\newcommand*{\NHonev}[2][\defaultdomain]{{\N{#2}}_{\Honev[{#1}]}}
\newcommand*{\SNHone}[2][\defaultdomain]{{\SN{#2}}_{\Hone[{#1}]}}

\newcommand*{\Hcurl}[1][\defaultdomain]{\boldsymbol{H}(\curl,{#1})}
\newcommand*{\bHcurl}[2][\defaultdomain]{\boldsymbol{H}_{#2}(\curl,{#1})}
\newcommand*{\zbHcurl}[1][\defaultdomain]{\bHcurl[#1]{0}}

\newcommand*{\NHcurl}[2][\defaultdomain]{\N{#2}_{\Hcurl[#1]}}

\newcommand{\D}{\mathrm{d}}

\newcommand{\ol}{\overline}

\newcommand{\be}{\begin{eqnarray}}
\newcommand{\ee}{\end{eqnarray}}

\newcommand{\ben}{\begin{eqnarray*}}
\newcommand{\een}{\end{eqnarray*}}

\newtheorem{theorem}{\sc Theorem}[section]
\newtheorem{remark}[theorem]{\sc Remark}
\newtheorem{algorithm}[theorem]{\sc Algorithm}
\newtheorem{example}[theorem]{\sc example}

\title{A robust solver for the finite element approximation of stationary incompressible MHD equations in 3D \thanks{}}

\author{Lingxiao Li
\thanks{Academy of Mathematics and System Sciences, Chinese Academy of Sciences; School of Mathematical Science, University of Chinese Academy of Sciences, Beijing, 100190, China.(lilingxiao@lsec.cc.ac.cn)}
\and Weiying Zheng
\thanks{Academy of Mathematics and System Sciences, Chinese Academy of Sciences; School of Mathematical Science, University of Chinese Academy of Sciences, Beijing, 100190, China.
This author is supported by China NSF under the grants 91430215 and by the National Magnetic
Confinement Fusion Science Program 2015GB110003.(zwy@lsec.cc.ac.cn)}}
\begin{document}
\date{}
\maketitle

\begin{abstract}
  In this paper, we propose a robust solver for the finite element discrete problem of the stationary incompressible magnetohydrodynamic (MHD) equations in three dimensions. By the mixed finite element method, both the velocity and the pressure are approximated by $\Honev$-conforming finite elements, while
  the magnetic field is approximated by $\Hcurl$-conforming edge elements.
  An efficient preconditioner is proposed to accelerate the convergence of the GMRES method for solving the linearized MHD problem. We use three numerical experiments to demonstrate the effectiveness of the finite element method and the robustness of the discrete solver. The preconditioner contains the least undetermined parameters and is optimal with respect to the number of degrees of freedom. We also show the scalability of the solver for moderate physical parameters.
\end{abstract}


\section{Introduction}
\label{sec:introduction}

Magnetohydrodynamics (MHD) has broad applications in our real world.
It describes the interaction between electrically conducting fluids
and magnetic fields. It is used in industry to heat, pump, stir, and
levitate liquid metals. Incompressible MHD model also governs the
terrestrial magnetic filed maintained by fluid motion in the earth
core and the solar magnetic field which generates sunspots and solar
flares\cite{dav01}. The incompressible MHD model consists of the
incompressible Navier-Stokes equations and the quasi-static Maxwell
equations. The magnetic field influences the momentum of the fluid
through Lorentz force, and conversely, the motion of fluid
influences the magnetic field through Faraday's law. In this paper,
we are studying the efficient iterative solver for the stationary MHD equations
\begin{subequations}\label{eq:mhd}
\begin{align}
 \Bu\cdot\nabla\Bu  +\nabla p -R_e^{-1}\Delta\Bu
 -S\curl\BB \times \BB
 =\Bf \qquad  \hbox{in}\;\;\Omega,   \label{eq:ns}\\
 \curl\left(\BB\times\Bu+R_m^{-1}\curl\BB\right)= 0
 \qquad \hbox{in}\;\;\Omega,  \label{eq:Faraday} \\
 \Div\Bu =0,   \quad\Div\BB=0 \qquad \hbox{in}\;\;\Omega,
\label{eq:div}
\end{align}
\end{subequations}
where $\Bu$ is the velocity of the fluid, $p$ is the hydrodynamic
pressure, $\BB$ is the magnetic flux density or the magnetic field provided with constant permeability, $R_e$ is the fluid Reynolds number,
$R_m$ is the magnetic Reynolds number, $S$ is the coupling
constant concerning the Lorentz force, and $\Bf\in\Ltwov$ stands for
the external force. We assume that $\Omega$ is a bounded Lipschitz domain. The system of equations are complemented with
Dirichlet boundary conditions
 \begin{equation}\label{eq:bc}
 \Bu=\Bg,\quad \BB\times\Bn =\BB_s\times\Bn \qquad\hbox{on}\;\;
 \Gamma:=\partial\Omega.
 \end{equation}

There are extensive papers in the literature to study numerical
solutions of incompressible MHD equations (cf. e.g.
\cite{ger06,gre10,gun91,hu15,la07,ma16,ni12,ni07-1,ni07-2,phi14}
and the references therein). In \cite{gun91}, Gunzburger et al
studied well-posedness and the finite element
method for the stationary incompressible MHD
equations. The magnetic field is discretized by the $\Honev$-conforming
finite element method. Strauss et al studied the adaptive finite
element method for two-dimensional MHD equations \cite{la07}.
Very recently, based on the nodal finite element approximation to
the magnetic field, Philips et al proposed a block preconditioner
based on an exact penalty formulation of stationary MHD equations
\cite{phi14}. We also refer to \cite{ger06} for a systematic
analysis on finite element methods for incompressible MHD
equations. When the domain has re-entrant angle, the magnetic field
may not be in $\Honev$. It is preferable to use noncontinuous finite
element functions to approximate $\BB$, namely, the so-called edge
element method\cite{hip02,ne86}. In 2004, Sch\"{o}tzau
proposed a mixed finite element method to solve the stationary
incompressible MHD equations where edge elements are used to solve the
magnetic field. To our knowledge, efficient solvers for
three-dimensional (3D) MHD equations are still rare in the literature,
particularly, for large Reynolds number $R_e$ and large coupling number $S$. An efficient solver should possess two merits:
\begin{enumerate}
\item the convergence rate is independent of the mesh or the number
of degrees of freedom (DOFs);

\item the algorithm is robust with respect to the physical parameters.
\end{enumerate}
The objective of this paper is to propose a preconditioned GMRES
method to solve the linearized discrete problem of \eqref{eq:mhd}--\eqref{eq:bc}. We shall adopt the mixed finite element method
proposed in \cite{sch04} and study efficient preconditioners for the
linearized problem.

Over the past three decades, fast solvers for incompressible
Navier-Stokes equations are relatively well-studied in the
literature (cf. e.g. \cite{elm06,elm07,elm14,ol99,zen94,zen95}).
For moderate Reynolds number, the Picard iteration for stationary
impressible Navier-Stokes equation is stable and efficient. At
each iteration, one needs to solve the linearized problem, the Oseen
equations
\begin{subequations}\label{eq:oseen}
\begin{align}
 \Bw\cdot\nabla\Bu  +\nabla p - R_e^{-1}\Delta\Bu =\Bf
 & \qquad  \hbox{in}\;\;\Omega,   \\
  \Div\Bu =0&   \qquad  \hbox{in}\;\;\Omega,   \\
 \Bu=\Bg& \qquad \hbox{on}\;\;\Gamma,
\end{align}
\end{subequations}
where $\Bw$ is the approximate solution at the previous step.
Iterative methods for discrete Oseen equations mainly consist of Krylov subspace methods,
multigrid methods, or their combinations. In terms of parallel computing and practical implementation, it is preferable to use
Krylov subspace method combined with an effective preconditioner. Among them, the pressure
convection-diffusion (PCD) preconditioner \cite{kay02}, the
least-squares commutator (LSC) preconditioner \cite{elm06,elm07,elm14}, and the augmented Lagrangian(AL) preconditioner \cite{be06,be11,ol03} prove robust and efficient for relatively large Reynolds number. In this paper, we shall
study the AL finite element method for the stationary MHD equations. Based on
a Picard-type linearization of the disrecte problem, we develop an efficient preconditioner for solving the linear problem. The preconditioner proves to be
robust when the Reynolds number and the coupling number are relatively large
and to be optimal with respect to the number of DOFs.

The paper is organized as follows:
In section 2, we introduce some notations for Sobolev spaces. A mixed finite element method is proposed to solve the AL formulation of the stationary MHD equations.
In section 3, we introduce a Picard-type linearization for the discrete MHD problem and devise an efficient preconditioner for solving the linear discrete problem. A preconditioned GMRES algorithm is also presented for the implementation of the discrete solver.
In Section 4, we present three numerical experiments to verify the optimal convergence rate of the mixed finite element method, to demonstrate the optimality and the robustness of the MHD solver, and to demonstrate the scalability for parallel computing.
Throughout the paper we denote vector-valued quantities by boldface notation, such as $\Ltwov:=(\Ltwo)^3$.

\section{Mixed finite element method for the MHD equations}
First we introduce some Hilbert spaces and Sobolev norms used in
this paper. Let $L^2(\Omega)$ be the usual Hilbert space of square
integrable functions equipped with the following inner product and
norm:
\begin{eqnarray}
(u,v):=\int_{\Omega}u(\Bx)\,v(\Bx)\D\Bx \quad \hbox{and} \quad
\NLtwo{u}:=(u,u)^{1/2}. \nonumber
\end{eqnarray}
Let the quotient space of $\Ltwo$ be defined by
 \ben
 L^2_0(\Omega):=\left\{ v\in\Ltwo\; : \; \int_\Omega v(\Bx)\D \Bx =0  \right\}
 = \Ltwo/\bbR\;.
 \een
Define $H^m(\Omega):=\{v\in L^2(\Omega): D^{\xi}v\in
L^2(\Omega),|\xi|\le m\}$ where $\xi$ represents non-negative triple
index. Let $H^1_0(\Omega)$ be the subspace of $H^1(\Omega)$ whose
functions have zero traces on $\Gamma$.

We define the spaces of functions having square integrable curl by
\begin{eqnarray}
\Hcurl&:=&\{\Bv\in\Ltwov\,:\;\curl\Bv\in \Ltwov\}, \nonumber\\
\zbHcurl&:=&\{\Bv\in\Hcurl\,:\;\Bn\times\Bv=0\;\;
\hbox{on}\;\Gamma\}, \nonumber
\end{eqnarray}
which are equipped with the following inner product and norm
\begin{equation*}
(\Bv,\Bw)_{\Hcurl}:=(\Bv,\Bw)+(\curl\Bv,\curl\Bw), \;\;
\NHcurl{\Bv}:=\sqrt{(\Bv,\Bv)_{\Hcurl}}\;.
\end{equation*}
Here $\Bn$ denotes the unit outer normal to $\Gamma$.

With a Lagrange multiplier $r$, we can rewrite \eqref{eq:mhd} into an AL form
\begin{subequations}\label{mhdr}
\begin{align}
 \Bu\cdot\nabla\Bu  +\nabla p - \gamma\nabla\Div\Bu
 - R_e^{-1}\Delta\Bu
 -S\curl\BB \times \BB
 =\Bf \qquad  \hbox{in}\;\;\Omega,   \label{up}\\
 S \curl\left(\BB\times\Bu + R_m^{-1}\curl\BB\right)
 +  \nabla r= 0
 \qquad \hbox{in}\;\;\Omega,  \label{Br} \\
 \Div\Bu =0,   \qquad \Div\BB=0 \qquad \hbox{in}\;\;\Omega,\label{div} \\
 \Bu =\Bg,   \quad\BB\times\Bn=\BB_s\times\Bn,\quad
 r=0 \qquad \hbox{on}\;\;\Gamma.\label{bc}
\end{align}
\end{subequations}
where $\gamma>0$ is the stabilization parameter or penalty parameter.
Taking divergence on both sides of \eqref{Br} and using \eqref{bc} yields
 \ben
 \Delta r =0\quad\hbox{in}\;\;\Omega,\qquad
 r=0\quad\hbox{in}\;\;\Gamma.
 \een
This means $r=0$ in $\Omega$ actually. Note that $\Div\Bu=0$, thus \eqref{mhdr} is equivalent to
\eqref{eq:mhd}--\eqref{eq:bc}. In the rest of this paper, we are
going to study the augmented problem \eqref{mhdr} instead of the
original problem.

A weak formulation of \eqref{mhdr} reads: Find
$(\Bu,\BB)\in\Honev\times\Hcurl$ and $(p,r)\in
L^2_0(\Omega)\times\zbHone$ such that $\Bu=\Bg$ and $\BB\times\Bn=\BB_s\times\Bn$ on $\Gamma$
and
\begin{subequations}\label{weak}
\begin{align}
\Ca((\Bu,\BB),(\Bv,\varphi)) + \Co((\Bu,\BB);(\Bu,\BB),(\Bv,\varphibf))
-\Cb((p,r),(\Bv,\varphibf)) &= (\Bf,\Bv),\label{weak:ab}\\
\Cb((q,s),(\Bu,\BB)) &= 0,\label{weak:b}
\end{align}
\end{subequations}
for all $(\Bv,\varphi)\in\zbHonev\times\zbHcurl$ and $(q,s)\in
L^2_0(\Omega)\times\zbHone$, where the bilinear forms and trilinear
form are defined respectively by
 \ben
 \Ca((\Bu,\BB),(\Bv,\varphi)) &=& R_e^{-1}(\nabla\Bu,\nabla\Bv)
    +\gamma(\Div\Bu, \Div\Bv)
    + SR_m^{-1}(\nabla\times\BB,\nabla\times\varphibf), \\
 \Co((\Bw,\psibf);(\Bu,\BB),(\Bv,\varphi))
  &=&(\Bw\cdot\nabla\Bu,\Bv) - S\left[(\curl\BB,\psibf\times\Bv)
    -(\psibf\times\Bu, \curl\varphibf)\right], \\
 \Cb((p,r),(\Bv,\varphibf)) &=& (p,\Div\Bv)
 + (\nabla r,\varphibf).
 \een
Assuming small data, Sch\"{o}tzau proved the existence and
uniqueness of the solution to \eqref{weak} without the penalized term $\gamma(\Div\Bu, \Div\Bv)$. The purpose of this paper is to propose a robust solver for the discrete problem.
This extra term in the new formula makes the discrete problem more well-defined for high Reynolds number\cite{ol03}.

Now we introduce the finite element approximation to \eqref{weak}. Let
$\Ct_h$ be a quasi-uniform and shape-regular tetrahedral mesh of
$\Omega$. Let $h$ denote the maximal diameter of all tetrahedra on
the mesh. For any $T\in\Ct_h$, let $P_k(T)$ be the space of
polynomials of degree $k\ge 0$ on $K$ and
$\BP_k(T)=\left(P_k(T)\right)^3$ be the corresponding space of
vector polynomials. Define the Lagrange finite element space of the
$k$-th order by
 \ben
 V(k,\Ct_h) =\left\{v\in\Hone:\;v|_T\in P_{k+1}(T),\;\forall\,T\in\Ct_h\right\}.
 \een
First we choose the well-known Taylor-Hood $P_2$-$P_1$ elements
\cite[Page 217-219]{bre91} for the discretization of $(\Bu,p)$,
namely,
 \ben
 \BV_h:= V(2,\Ct_h)^3\cap\zbHonev,\qquad Q_h:=V(1,\Ct_h)\;.
 \een
From \cite[Page 255-258]{bre91}, the discrete inf-sup condition
holds
 \begin{equation}\label{infsup-u}
 \sup_{0\ne \Bv\in\BV_h} \frac{(q,\Div\Bv)}{\NHone{\Bv}} \ge C_u \NLtwo{q}
 \qquad\forall\,q\in Q_h,
 \end{equation}
where $C_u$ is the inf-sup constant independent of the mesh
size. We shall also use $ \ol\BV_h = V(k,\Ct_h)^3$.

The finite element space for $\BB$ is chosen as N\'{e}d\'{e}lec's
edge element space of the first order in the second family\cite{ne86}, namely,
 \ben
 \ol\BC_h =\left\{\Bv\in\Hcurl:\;\Bv|_T\in\BP_1(T),
 \;\forall\,T\in\Ct_h\right\},\qquad
 \BC_h=\ol\BC_h\cap\zbHcurl.
 \een
The finite element space for $r$ is defined by
 \ben
 S_h =V(2,\Ct_h)\cap\zbHone.
 \een
Since $\nabla S_h\subset\BC_h$, we easily get the inf-sup condition for the pair of finite element spaces $\BC_h\times S_h$
 \be\label{infsup-B}
  \sup_{0\ne \Bv\in\BC_h} \frac{(\nabla s,\Bv)}{\NHcurl{\Bv}} \ge \SNHone{s}\ge
  C_b \NHone{s}\qquad\forall\,s\in S_h,
 \ee
where $C_b>0$ is the Poinc\'{a}re constant depending only on $\Omega$.

The finite element approximation to \eqref{weak} reads: Find
$(\Bu_h,\BB_h)\in\ol\BV_h\times\ol\BC_h$ and $(p_h,r_h)\in Q_h\times S_h$
such that
\begin{subequations}\label{weakh}
\begin{align}
 \Ca((\Bu_h,\BB_h),(\Bv,\varphibf)) + \Co((\Bu_h,\BB_h);(\Bu_h,\BB_h),(\Bv,\varphibf))
  -\Cb((p_h,r_h),(\Bv,\varphibf)) &= (\Bf,\Bv),
  \label{weakh:ab}\\
  \Cb((q,s),(\Bu_h,\BB_h)) &= 0,\label{weakh:b}
\end{align}
\end{subequations}
for all $(\Bv,\varphibf)\in\BV_h\times\BC_h$ and $(q,s)\in Q_h\times
S_h$. From \eqref{infsup-u} and \eqref{infsup-B} we know that the
bilinear form $\Cb(\cdot,\cdot)$ satisfies the discrete inf-sup
condition
 \be\label{infsup-h}
  \sup_{(\Bv_h,\varphibf_h)\in\BV_h\times\BC_h}
  \frac{\Cb((q_h,s_h),(\Bv_h,\varphibf_h))}{\N{(\Bv_h,\varphibf_h)}_{\BV_h\times\BC_h}} \ge\min(C_u,C_p)\N{(q_h,s_h)}_{Q_h\times S_h}
  \quad\forall\,(q_h,s_h)\in Q_h\times S_h\;,
 \ee
where
 \ben
 \N{(\Bv_h,\varphibf_h)}_{\BV_h\times\BC_h}
    :=\sqrt{\NHonev{\Bv_h}^2+\NHcurl{\varphibf_h}^2}\;,\qquad
 \N{(q_h,s_h)}_{Q_h\times S_h}:=\sqrt{\NLtwo{q_h}^2+\NHone{s_h}^2}\;.
 \een
Based on the assumption of small data, we can prove that the discrete problem \eqref{weakh} has a unique solution.
Again we do not elaborate on the details and pay our attention to fast solvers of the discrete solution \eqref{weakh}.


\section{A preconditioner for the linearized finite element problem}

In this section, we are going to study the solution of the nonlinear
discrete problem \eqref{weakh}. First we propose a Picard-type
iterative method for solving \eqref{weakh}. At each nonlinear iteration, the linearized problem consists of an AL Oseen equation with Lorentz force and a Maxwell equation coupled with the fluid. The preconditioner for the
linearized MHD equation depends crucially on the preconditioner for
the penalized Navier-Stokes equations and the preconditioner for the Maxwell
equations in mixed forms.

\subsection{Picard-type method for the discrete MHD equations}

In this subsection, we consider the Picard linearization of \eqref{weakh}.
For convenience, we rearrange the order of variables as $(\BB_h, r_h, \Bu_h, p_h)$ in the linearized problem.
Let
$(\BB_{k},r_{k}, \Bu_{k}, p_{k})\in \ol\BC_h\times S_h \times
\ol\BV_h\times Q_h$
be the approximate solutions of \eqref{weakh} from the previous iteration. The error equation for these approximate solutions reads:
Find $(\delta \BB_{k}, \delta r_{k}, \delta \Bu_k, \delta p_{k})\in\BV_h\times Q_h\times\BC_h\times S_h$ such that
\begin{subequations}\label{Linear}
\begin{align}
 SR_m^{-1}(\curl\delta\BB_{k}, \curl\varphibf) + (\nabla\delta r_{k}, \varphibf) + S(\BB_{k}\times\delta\Bu_k, \curl\varphibf) = R_b(\varphibf), &\label{Linear:b} \\
 (\delta \BB_{k}, \nabla s) = R_r(s), &\label{Linear:r}\\
 -S(\curl\delta\BB_{k}, \BB_k\times\Bv) + \mathcal{F}(\Bu_k;\delta\Bu_k,\Bv)- (\delta p_{k}, \Div\Bv) = R_u(\Bv), &\label{Linear:u}\\
 -(\Div\delta \Bu_k, q) = R_p(q), &\label{Linear:p}
\end{align}
\end{subequations}
where the trilinear form $\Cf$ represents the convection-diffusion part of the fluid equation
 \ben
 \mathcal{F}(\Bu_k;\delta\Bu_k,\Bv) := R_e^{-1}(\nabla\delta \Bu_k, \nabla\Bv)
 + (\Bu_k\cdot\nabla\delta\Bu_k,\Bv) +\gamma(\Div\delta \Bu_k, \Div\Bv),
 \een
and the residual functionals are defined by
\begin{align*}
 &R_b(\varphibf) = -SR_m^{-1}(\curl\BB_{k}, \curl\varphibf) - (\nabla r_{k}, \varphibf) - S(\BB_{k}\times\Bu_k, \curl\varphibf),\\
 &R_r(s) = - (\BB_k, \nabla s),\\
 &R_u(\Bv) = (\Bf,\Bv) - \mathcal{F}(\Bu_k;\Bu_k,\Bv) + S(\curl\BB_{k}, \BB_k\times\Bv) + (p_k, \Div\Bv),\\
 &R_p(q) = (\Div\Bu_k, q).
\end{align*}
After solving \eqref{Linear}, the approximate solutions will be
updated by
\begin{equation}
\BB_{k+1}=\BB_{k}+\theta\delta\BB_{k},~~ r_{k+1}=r_{k}+\theta\delta r_{k},~~ \Bu_{k+1}=\Bu_k + \theta\delta \Bu_k,~~ p_{k+1}=p_k + \theta\delta p_{k}
\end{equation}
with a relaxation factor $0<\theta\le 1$.

To devise the preconditioner, we write problem \eqref{Linear} into an algebraic form
\be\label{Axb}
 \bbA\mathbf{x} = \mathbf{b}\;,
\ee
where the solution vector $\mathbf{x}$ consists of the degrees of
freedom for $(\delta \BB_{k}, \delta r_{k}, \delta \Bu_k, \delta p_{k})$ respectively,
$\mathbf{b}$ is the residual vector, and $\bbA$ is the stiffness matrix.
In block forms, they can be written as
\be\label{matrixA}
 \mathbf{x}=
 \left(
 \begin{array}{c}
 \mathbf{x}_b \\  \mathbf{x}_r \\  \mathbf{x}_u \\  \mathbf{x}_p
 \end{array}
 \right),\qquad
 \mathbf{b}=
 \left(
 \begin{array}{c}
 R_b \\  R_r \\  R_u \\  R_p
 \end{array}
 \right),\qquad
 \bbA=
 \left(
 \begin{array}{cccc}
 \bbC  &\bbG^\top  &\bbJ^\top  &0 \\
 \bbG  &0          &0          &0 \\
 -\bbJ &0          &\bbF       &\bbB^\top \\
 0     &0          &\bbB       &0
 \end{array}
 \right).
\ee
Let $\{\Bv_i: \;1\le i\le N_V\}$,  $\{\varphibf_i: \;1\le i\le
N_C\}$, $\{q_i: \;1\le i\le N_Q\}$,  $\{s_i: \;1\le i\le N_S\}$ be
the bases of $\BV_h$, $\BC_h$, $Q_h$, and $S_h$ respectively. Then
the entries of all block matrices are defined by
\begin{align*}
 \bbC_{ij} &= SR_m^{-1}(\curl\varphibf_j,\curl\varphibf_i)\;, \\
 \bbG_{ij} &= S(\varphibf_j, \nabla s_i)\;,\\
 \bbJ_{ij} &= S(\curl\varphibf_j,\BB_k\times\Bv_i)\;\\
 \bbF_{ij} &= \mathcal{F}(\Bu_k;\Bv_j, \Bv_i)\;,\\
 \bbB_{ij} &= -(\Div\Bv_j, q_i)\;.
\end{align*}
Clearly the block matrices represents the differential operators appearing in the Navier-Stokes equations and the Maxwell equation
on various finite element spaces
\begin{align*}
 \bbC\Leftrightarrow SR_m^{-1}\curl\curl, \quad
  \bbG \Leftrightarrow -\Div
  & \qquad\hbox{on}\;\;\BC_h, \\
    \bbG^\top \Leftrightarrow \nabla
  & \qquad\hbox{on}\;\;S_h, \\
 \bbF\Leftrightarrow \left(-R_e^{-1}\Delta+\Bu_k\cdot\nabla
 	-\gamma\nabla\Div\right), \quad
   \bbB \Leftrightarrow -\Div
 &\qquad\hbox{on}\;\;\BV_h,\\
    \bbB^\top \Leftrightarrow \nabla
  & \qquad\hbox{on}\;\;Q_h.
\end{align*}
Here $-\Div$ is understood as the dual operator of $\nabla|_{S_h}$ or $\nabla|_{Q_h}$.
Moreover, $\bbJ, \bbJ^\top$ are algebraic representations of
the two multiplication operators which couple the magnetic field and the conducting fluid. For any given $\Bw\in\BV_h$ and $\psibf\in\BC_h$, we have
 \be\label{oper-JJT}
 \bbJ^\top \Leftrightarrow S\curl(\BB_k\times \Bw)
 			\quad \hbox{on}\;\;\BC_h, \qquad
 \bbJ \Leftrightarrow S\curl\psibf\times\BB_k
 			\quad \hbox{on}\;\;\BV_h.
 \ee
The relationships between these operators play an important role in
devising a robust preconditioner for the linearized problem.

\subsection{Preconditioning for the linearized MHD equations}

Let $\bbL_r$ be the stiffness matrix of $-\Delta$ on $S_h$ and let $\sigma>0$ be a constant. First we post-multiply the second column of $\bbA$ by $\sigma\bbL_r^{-1}\bbG$ and add it to the first column. This yields a matrix
\ben
 \bbA_1= \left(
 \begin{array}{cccc}
 \bbC + \sigma\bbS_r  & \bbG^\top  & \bbJ^\top &  0 \\
 \bbG  &0 &0     &0 \\
 -\bbJ &0 &\bbF  &\bbB^\top \\
 0     &0 &\bbB  &0
 \end{array}
 \right), \qquad
 \bbS_r := \bbG^\top\bbL_r^{-1}\bbG.
\een
Next, pre-multiplying the first row of $\bbA_1$ by $-\bbG(\bbC + \sigma\bbS_r)^{-1}$ and adding it to the second row, we get a matrix
\ben
 \bbA_2= \left(
 \begin{array}{cccc}
 \bbC + \sigma\bbS_r              & \bbG^\top  & \bbJ^\top      &0 \\
 0     &-\bbG(\bbC + \sigma \bbS_r)^{-1}\bbG^\top &-\bbG(\bbC + \sigma\bbS_r)^{-1}\bbJ^\top  &0 \\
 -\bbJ &0 &\bbF  &\bbB^\top \\
 0     &0 &\bbB  &0
 \end{array}
 \right).
\een
Note that $\bbG(\bbC+\sigma \bbS_r)$ represents the operator
$\Div\left(SR_m^{-1}\curl\curl+\sigma\nabla\Delta^{-1}\Div\right)$. Since
 \ben
 \Div \left(SR_m^{-1}\curl\curl+\sigma\nabla\Delta^{-1}\Div\right)
 =\sigma\Delta \Delta^{-1}\Div =\sigma\Div,
 \een
formally we have
 \ben
 \Div\left(SR_m^{-1}\curl\curl+\sigma\nabla\Delta^{-1}\Div\right)^{-1}
  =\sigma^{-1}\Div .
 \een
This means that
 \ben
 \bbG(\bbC+\sigma \bbS_r)^{-1} \approx \sigma^{-1}\bbG .
 \een
Since $\bbJ^\top$ represents the coupling term $\curl(S\BB_k\times\Bv_h)$, we have $\bbG\bbJ^\top\approx 0$. Therefore,
 \ben
 \bbG(\bbC+\sigma \bbS_r)^{-1}\bbJ^\top \approx 0.
 \een
Moreover, from \cite{gre07}, we know
that $\bbC + \sigma\bbM$ is equivalent to $\bbC + \sigma\bbG^\top\bbL_r^{-1}\bbG$ in spectrum where $\bbM$ is the mass matrix on $\BC_h$. So one gets the approximation
\ben
 \bbA_2\approx \bbA_3:=  \left(
 \begin{array}{cccc}
 \bbC + \sigma\bbM  & \bbG^\top  & \bbJ^\top &  0 \\
 0 & -\sigma^{-1}\bbL_r & 0 & 0 \\
 -\bbJ & 0 &\bbF  & \bbB^\top \\
 0 & 0 & \bbB & 0
 \end{array}
 \right) .
\een

Next, pre-multiplying the first row of $\bbA_3$ by $\bbJ(\bbC + \sigma\bbM)^{-1}$ and adding it to the third row, we get a matrix
\ben
 \bbA_4 = \left(
 \begin{array}{cccc}
 \bbC + \sigma\bbM  & \bbG^\top  & \bbJ^\top &  0 \\
 0 & -\sigma^{-1}\bbL_r & 0 & 0 \\
 0 & \bbJ(\bbC + \sigma\bbM)^{-1}\bbG^\top &\bbF + \bbJ(\bbC + \sigma\bbM)^{-1}\bbJ^\top  & \bbB^\top \\
 0 & 0 & \bbB & 0
 \end{array}
 \right).
\een
Since $\Div\left(SR_m^{-1}\curl\curl +\sigma\VI\right)=\sigma\Div$, formally we have
$\Div \left(SR_m^{-1}\curl\curl +\sigma\VI\right)^{-1} =\sigma^{-1}\Div$. Here $\VI$ is the identity operator. This means
 \ben
 \bbG (\bbC+\sigma\bbM)^{-1}\bbJ^\top \approx 0
 \qquad \hbox{or}\qquad
 \bbJ (\bbC+\sigma\bbM)^{-1}\bbG^\top \approx 0.
 \een
So one obtains an approximation of $\bbA_4$ as follows
 \be\label{A5}
 \bbA_4 \approx \bbA_5:= \left(
 \begin{array}{cccc}
 \bbC + \sigma\bbM  & \bbG^\top  & \bbJ^\top &  0 \\
 0 & -\sigma^{-1}\bbL_r & 0 & 0 \\
 0 & 0 &\bbF + \bbJ(\bbC + \sigma\bbM)^{-1}\bbJ^\top  & \bbB^\top \\
 0 & 0 & \bbB & 0
 \end{array}
 \right).
 \ee
This means that $\bbA_5^{-1}$ is actually a natural preconditioner for $\bbA$,
except for the difficulties in computing the inverse of the block $\bbF + \bbJ(\bbC + \sigma\bbM)^{-1}\bbJ^\top$. In the next subsection, we are going to
derive a good approximation of
$\bbF + \bbJ(\bbC + \sigma\bbM)^{-1}\bbJ^\top$ so that its approximation inverse is easy to compute iteratively.

\subsection{An efficient preconditioner for the magnetic field-fluid coupling block}
\label{subsec:coupling}

From \eqref{A5}, the key step to compute $\bbA_5^{-1}$
is how to precondition the $2\times 2$ block
 \be\label{matrixX}
 \bbX=\left(
 \begin{array}{cc}
 \bbC + \sigma\bbM   & \bbJ^\top  \\
  0  &\bbF + \bbJ(\bbC+\sigma\bbM)^{-1}\bbJ^\top
 \end{array}
 \right).
 \ee
We note that $\bbF + \bbJ(\bbC + \sigma\bbM)^{-1}\bbJ^\top$ is the precise Schur complement of the following matrix which accounts for the coupling between $\BB_h$ and $\Bu_h$
 \be\label{hatX}
 \hat\bbX=\left(
 \begin{array}{cc}
 \bbC + \sigma\bbM   & \bbJ^\top  \\
  -\bbJ  &\bbF
 \end{array}
 \right).
 \ee
Remember from \eqref{oper-JJT} that $\bbJ^\top$ and $\bbJ$
represent, respectively, the two multiplication operators
$S\curl(\BB_k\times \Bw)$ and $S\curl\psibf\times\BB_k$ for any given $\Bw$ and $\psibf$.
So $\bbJ (\bbC + \sigma\bbM)^{-1}\bbJ^\top\Bw$ is the algebraic representation of
\begin{equation}\label{BuSchur}
S\curl\left\{\left(SR_m^{-1}\curl\curl
	+ \sigma \mathbf{I}\right)^{-1}S\curl(\BB_k\times \Bw)
	\right\}\times\BB_k
\end{equation}
Since $(S R_m^{-1}\curl\curl + \sigma \mathbf{I})^{-1}$  commutates with $\curl$, then \eqref{BuSchur} becomes
\ben
S^2 \left\{\left(SR_m^{-1}\curl\curl + \sigma \mathbf{I}\right)^{-1}
 \curl\curl(\BB_k\times \Bw)\right\}\times\BB_k
\een
For $\sigma>0$ sufficiently small, we adopt the following approximation
 \ben
  \curl\curl  \approx  S^{-1}R_m(SR^{-1}_m\curl\curl + \sigma\mathbf{I}) .
  \een
This yields an approximation of \eqref{BuSchur}
\begin{equation}\label{JinsiSchur}
 S^2\left\{\left(SR^{-1}_m\curl\curl
	+ \sigma \mathbf{I}\right)^{-1}\curl\curl(\BB_k\times \Bw)
	\right\}\times\BB_k \approx
 SR_m(\BB_k\times \Bw)\times\BB_k .
\end{equation}
Therefore, we get an approximation of the Schur complement
 \ben
 \bbF + \bbJ(\bbC+\sigma\bbM)^{-1}\bbJ^\top
 \approx \bbS,
 \een
where $\bbS$ is the stiffness matrix associated with the bilinear
form
 \ben
  \mathcal{F}(\Bu_k;\Bu,\Bv) +
  SR_m(\BB_k\times \Bu,\BB_k\times\Bv) .
 \een
This should yield a good preconditioner of $\bbX$, that is,
\begin{equation}\label{preub}
  \left(
 \begin{array}{cc}
 \bbC + \sigma\bbM   & \bbJ^\top  \\
  0  &\bbF + \bbJ(\bbC+\sigma\bbM)^{-1}\bbJ^\top
 \end{array}
 \right) \sim
 \left(
 \begin{array}{cc}
 \bbC + \sigma\bbM   & \bbJ^\top  \\
 0  &\bbS
 \end{array}
 \right).
\end{equation}

\begin{remark}
In \cite{phi14}, Philips et al studied the preconditioner for two-dimensional MHD equations where the magnetic field is discretized with $\Hone$-conforming finite elements.
Based on an exact penalty formulation, they propose to approximate the coupling effect between magnetic field and velocity by
 \ben
 \beta SR_m(\BB_k\times \Bw)\times\BB_k,
 \een
where $\beta>0$ is a parameter depending on
the mesh size and the magnitude of $\BB_k$.

Compared with \cite{phi14}, for the 3D MHD equations and the $\Hcurl$-conforming approximation of $\BB$, our approximation to the coupling effect in the preconditioner level
does not need the extra parameter and has more advantages in practical computations such as adaptive computing, though they are similar.
\end{remark}

Now we demonstrate numerically the robustness of the preconditioner in \eqref{preub} with respect to the parameter $\sigma$ and the mesh size $h$.  Since it is the approximation $\bbJ(\bbC+\sigma\bbM)^{-1}\bbJ^\top
\approx \bbS$ that is concerned here, we fix $R_e = 1.0$ and $\gamma = 1.2$ and test the efficiency of the preconditioner for
different values of $S$ and $R_m$.

We consider the system of linear equations: Find $\delta\Bu\in\BV_h$ and $\delta\BB\in\BC_h$ such that
\begin{subequations}\label{testub}
\begin{align}
-S(\curl\delta\BB, \BB_0\times\Bv) + \mathcal{F}(\Bu_0;\delta\Bu,\Bv) = (\Bf, \Bv)
 &\quad \forall\,\Bv\in\BV_h,  \\
 SR^{-1}_m(\curl\delta\BB, \curl\varphibf)
   + \sigma(\delta\BB, \varphibf)
   + S(\BB_{0}\times\delta\Bu_k, \curl\varphibf) = 0
  &\quad \forall\,\varphibf\in\BC_h,
\end{align}
\end{subequations}
where
 \ben
  \Bf = (1,\sin(x),0),\quad
  \Bu_0 = (y, \sin(x+z), 1), \quad
  \BB_0=(\sin(y)+\cos(z), 1-\sin(x), 1).
 \een
Clearly the stiffness matrix of \eqref{testub} is $\hat\bbX$ which is given in \eqref{hatX}. In the following, we test three cases of $\sigma$
 \ben
 \sigma = 1,\;  10^{-2},\; 10^{-4},
 \een
and three cases of physical parameters
 \ben
 S = R_m  = 1,\; 10,\; 100.
 \een
The computational domain is the unit cube, namely,
$\Omega=(0,1)^3$.

\begin{table}[!htb]
  \centering
  \caption{Number of preconditioned GMRES iterations for $\sigma=1$.}
  \label{test-table1}
  \begin{tabular*}{10cm}{@{\extracolsep{\fill}}|c|c|c|c|}
  \hline
  $h$  &$S=R_m=1$  &$S=R_m=10$ &$S=R_m=100$  \\
  \hline
  0.216506     &4  &14   &91   \\    \hline
  0.108253     &4  &14   &74   \\    \hline
  0.054127     &4  &14   &64   \\    \hline
  0.027063     &4  &14   &61   \\    \hline
  \end{tabular*}
\end{table}
\begin{table}[!htb]
  \begin{center}
  \caption{Number of preconditioned GMRES iterations for $\sigma=10^{-2}$.}\label{test-table3}
  \begin{tabular*}{10cm}{@{\extracolsep{\fill}}|c|c|c|c|}
  \hline
  $h$    &$S=R_m=1$  &$S=R_m=10$    &$S=R_m=100$
    \\  \hline
  0.216506     &4  &14   &91   \\    \hline
  0.108253     &4  &14   &73   \\    \hline
  0.054127     &4  &14   &64   \\    \hline
  0.027063     &4  &14   &61   \\    \hline
  \end{tabular*}
  \end{center}
\end{table}

\begin{table}[!htb]
  \begin{center}
  \caption{Number of preconditioned GMRES iterations for $\sigma=10^{-4}$.}\label{test-table4}
  \begin{tabular*}{10cm}{@{\extracolsep{\fill}}|c|c|c|c|}
  \hline
  $h$    &$S=R_m=1$  &$S=R_m=10$    &$S=R_m=100$
  \\  \hline
  0.216506     &4  &14   &91   \\    \hline
  0.108253     &4  &14   &73   \\    \hline
  0.054127     &4  &14   &64   \\    \hline
  0.027063     &4  &14   &61   \\    \hline
  \end{tabular*}
  \end{center}
\end{table}

\begin{table}[!htb]
  \begin{center}
  \caption{Number of preconditioned GMRES iterations for
  $\sigma=10^{-4}$ with $\bbS=\bbF$.}
  \label{test-table5}
  \begin{tabular*}{10cm}{@{\extracolsep{\fill}}|c|c|c|c|}
  \hline
  $h$    &$S=R_m=1$  &$S=R_m=10$    &$S=R_m=100$  \\  \hline
  0.216506     &3  &14   &82   \\    \hline
  0.108253     &4  &14   &92   \\    \hline
  0.054127     &4  &15   &97   \\    \hline
  0.027063     &4  &15   &98   \\    \hline
  \end{tabular*}
  \end{center}
\end{table}

We use preconditioned GMRES method to solve \eqref{testub} and the preconditioner is set by \eqref{preub}.
This means that we need solve the residual equation at each GMRES iteration
 \begin{equation}\label{resbu}
 \bbS\Ve_u =\Vr_u,\qquad
 (\bbC+\sigma\bbM)\Ve_b =\Vr_b -\bbJ^\top\Ve_u,
 \end{equation}
where $\Vr_b,\Vr_u$ stand for the residual vectors and
$\Ve_b,\Ve_u$ stand for the error vectors. The tolerance for
the relative residual of the GMRES method is set by $10^{-6}$.
The tolerances for solving the two sub-problems in \eqref{resbu} are set by $10^{-3}$.
From Table~\ref{test-table1}--\ref{test-table4}, we find that the convergence of the preconditioned GMRES is uniform with respect to both $\sigma$ and $h$. An interesting observation is that, for large $S=R_m$,  the number of GMRES iterations even decreases when $h\to 0$. In this case, the magnetic field-fluid
coupling becomes strong.

Table~\ref{test-table5} shows the number of preconditioned GMRES iterations where the approximate Schur complement $\bbS$ in \eqref{preub} is replaced with the matrix $\bbF$. It amounts to devise a preconditioner of $\bbX$ by dropping its left lower block
 \be\label{preXF}
 \left(
 \begin{array}{cc}
 \bbC + \sigma\bbM   & \bbJ^\top  \\
  -\bbJ  &\bbF
 \end{array}
 \right) \sim \left(
 \begin{array}{cc}
 \bbC + \sigma\bbM   & \bbJ^\top  \\
 0  &\bbF
 \end{array}
 \right) .
 \ee
This is the classical Riesz map preconditioning in \cite{ki10} or the operator preconditioning in \cite{hip06}. Comparing Table~\ref{test-table4} with Table~\ref{test-table5},
we find that, for large $S=R_m$, the convergence of GMRES method with this preconditioner becomes slower and deteriorates  when $h\to 0$. This becomes even more apparent when solving the whole MHD system (see Table~\ref{cavity:pre} for the computation of driven cavity flow).

\subsection{A preconditioner for the augmented Navier-Stokes equations}

Combining \eqref{A5} and \eqref{preub}, we get a preconditioner of $\bbA$, that is, the inverse of
 \be\label{preA}
 \bbA_6 := \left(
 \begin{array}{cccc}
 \bbC + \sigma\bbM  & \bbG^\top  & \bbJ^\top &  0 \\
 0 & -\sigma^{-1}\bbL_r & 0 & 0 \\
 0 & 0 &\bbS  & \bbB^\top \\
 0 & 0 & \bbB & 0
 \end{array}
 \right) .
 \ee
It is left to study the preconditioner for the lower right $2\times 2$ block of $\bbA_6$, namely,
\begin{align*}
 \left(
 \begin{array}{cc}
 \bbS & \bbB^\top  \\
 \bbB & 0
 \end{array}
 \right).
\end{align*}
It amounts to solve the saddle point problem:
Find $(\delta\Bu,\delta p)\in\BV_h\times Q_h$ such that
\begin{subequations}\label{eq:diff-oseen}
\begin{align}
 \mathcal{F}(\Bu_k;\delta\Bu,\Bv)
 +SR_m (\BB_k\times\delta\Bu,\BB_k\times\Bv)
 - (\delta p,\Div\Bv) &= R_u(\Bv) \qquad  \forall\,\Bv\in\BV_h,   \\
  -(\Div\delta\Bu,q) &=R_p(q)    \qquad  \forall\, q\in Q_h.
\end{align}
\end{subequations}
In \cite{be06, be11}, Benzi et al studied the Osceen equation(namely $\BB_k = \mathbf{0}$) and proposed to use the following preconditioner
\begin{equation}\label{NSschur}
 \left(
  \begin{array}{cc}
    \bbF & \bbB^\top \\
    0 &  -\left(R_e^{-1}+\gamma\right)^{-1}\bbQ_p\\
  \end{array}
\right)^{-1},
\end{equation}
where $\bbQ_p$ is the mass matrix on $Q_h$. It is proved that
the above preconditioner is efficient for relatively large Reynolds number. With $\Bu_k = \BB_k = \mathbf{0}$, namely the Stokes equations, we refer to \cite{lo04, ma11} for similar arguments.
And with $\Bu_k=\mathbf{0}$ for the time-dependent incompressible MHD, the work in \cite{ma16} give useful insight using pressure mass matrix as a subblock.
Inspired by them, we propose to precondition
\begin{align}\label{MNS}
 \left(
 \begin{array}{cc}
 \bbS & \bbB^\top  \\
 \bbB & 0
 \end{array}
 \right) \quad \hbox{by} \quad \left(
  \begin{array}{cc}
    \bbS & \bbB^\top \\
    0 &  -\left(R_e^{-1}+\gamma\right)^{-1}\bbQ_p\\
  \end{array}
\right)^{-1}.
\end{align}

\subsection{A robust preconditioner for the linearized MHD problem}

Using \eqref{preA} and \eqref{MNS}, a preconditioner for $\bbA$ is given by the inverse of
\begin{equation}\label{preII}
 \left(
 \begin{array}{cccc}
 \bbC + \sigma\bbM  & \bbG^\top  & \bbJ^\top &  0 \\
 0 & -\sigma^{-1}\bbL_r & 0 & 0 \\
 0 & 0 &\bbS  & \bbB^\top \\
 0 & 0 & 0    & -\left(R_e^{-1}+\gamma\right)^{-1}\bbQ_p
 \end{array}
 \right),
\end{equation}
where the parameter $\gamma$ can be used to tune the efficiency of the preconditioner. According to our experience,
any value $\gamma \sim O(1)$ works well for high Reynolds and moderate $S$ and $R_m$.
Moreover, to fix the parameter $\sigma$, we set $\sigma=SR^{-1}_m$ so that $\bbC+\sigma\bbM$ is associated with the bilinear form
 \ben
 SR^{-1}_m \left[(\curl\Bu,\curl\Bv) + (\Bu,\Bv) \right].
 \een
This yields our final preconditioner for the stiffness matrix
$\bbA$, defined by
\begin{equation}\label{preconditioner}
 \bbP= \left(
 \begin{array}{cccc}
 \bbC + SR^{-1}_m\bbM  & \bbG^\top  & \bbJ^\top &  0 \\
 0 & -S^{-1}R_m\bbL_r & 0 & 0 \\
 0 & 0 &\bbS  & \bbB^\top \\
 0 & 0 & 0 & -\left(R_e^{-1}+\gamma\right)^{-1}\bbQ_p
 \end{array}
 \right)^{-1}.
\end{equation}

Now we are in the position to present the preconditioned GMRES algorithm for solving the linear system \eqref{Axb} of the MHD problem. The idea is to use an approximation of $\bbP$ to precondition $\bbA$.
For convenience in notation, given a vector $\Vx$ which has the same size as one column vector of $\bbA$, we let $(\Vx_b, \Vx_r, \Vx_u, \Vx_p)$ be the vectors which consist of entries of $\Vx$ and correspond to $(\BB_h, r_h, \Bu_h, p_h)$ respectively.

\begin{algorithm} [Preconditioned GMRES Algorithm]\label{alg:gmres}
{\sf
Given the tolerances $\varepsilon\in
(0,1)$ and $\varepsilon_0\in (\varepsilon,1)$, the maximal number of GMRES iterations $N>0$, and the
initial guess $\mathbf{x}^{(0)}$ for the solution of \eqref{Axb}.
Set $k=0$ and compute the residual vector
\begin{align*}
\mathbf{r}^{(k)} = \mathbf{b}-\bbA\mathbf{x}^{(k)}.
\end{align*}

\noindent
While $\left(k<N\;\; \&\;\; \N{\mathbf{r}^{(k)}}_2
	>\varepsilon\N{\mathbf{r}^{(0)}}_2\right)$ do
\begin{enumerate}
\item Solve $\bbQ_p \Ve_ p = -\left(R_e^{-1}+\gamma\right)\Vr_p^{(k)}$ by the CG method with the diagonal preconditioning and tolerance $\varepsilon_0$
for the relative residual.
\vspace{1mm}

\item Solve $\bbS \Ve_u = \Vr^{(k)}_u - \bbB^\top \Ve_p$ by preconditioned GMRES method with tolerance $10^{-3}$.
The preconditioner is the one level additive Schwarz method with overlap = 2 \cite{Cai99}.
\vspace{1mm}

\item Solve $\bbL_r \Ve_r = -SR^{-1}_m\Vr^{(k)}_r$ by preconditioned CG method with tolerance $\varepsilon_0$. The preconditioner is the algebraic multigrid method (AMG) solver \cite{hen02}.
\vspace{1mm}

\item Solve $(\bbC + SR^{-1}_m\bbM)\Ve_b
= \Vr^{(k)}_b - \bbJ^\top \Ve_u - \bbG^\top \Ve_r$ by preconditioned CG method with tolerance $\varepsilon_0$ and the Hiptmair-Xu preconditioner \cite{hip07}.
\vspace{1mm}

\item Update the solution:
$\mathbf{x}^{(k+1)}:=\mathbf{x}^{(k)} + \mathbf{e}^{(k)}$ .
\vspace{1mm}

\item Set $k:=k+1$ and compute the residual vector
$\mathbf{r}^{(k)} = \mathbf{b}-\bbA\mathbf{x}^{(k)}$.
\end{enumerate}
End while.}
\end{algorithm}


\section{Numerical experiments}

In this section, we present three numerical experiments to verify the convergence rate of finite element approximation to the augmented Lagrangian(AL) formulation of the MHD,
to demonstrate the robustness of the preconditioner, and to demonstrate the scalability of the parallel solver.
The parallel code is developed based on the finite element package---Parallel Hierarchical Grids (PHG) \cite{zh09-1,zh09-2}.

\begin{example}
This example is to verify the convergence rate of the finite element
discrete problem \eqref{weakh}.
The analytic solutions are chosen as
 \ben
 \Bu=\left(
        \begin{array}{c}
          \sin z \\
          2\cos x \\
          0 \\
        \end{array}
      \right),\quad
p=\sin y + \cos 1-1,\quad
\BB=\left(
      \begin{array}{c}
        \cos y \\
        0 \\
        0 \\
      \end{array}
    \right),\quad
r=0.
 \een
The parameters are set by $R_e = S = R_m = 1$ and $\gamma = 1$.
\end{example}

\begin{table}[!htb]
  \centering
  \caption{Convergence rate of the finite element discrete problem.}
  \label{Error}
  \begin{tabular*}{14cm}{@{\extracolsep{\fill}}|c|c|c|c|c|c|c|}
  \hline
  $h$        &$\|\Bu-\Bu_h\|_{H^1}$   &order &$\|p - p_h\|_{L^2}$  &order &$\|\BB-\BB_h\|_{\BH(\curl)}$ &order\\ \hline
  0.4330   &2.893e-03   &---&1.848e-03   &---&4.811e-02  &---\\    \hline
  0.2165   &7.071e-04   &2.033&3.908e-04   &2.242&2.378e-02  &1.017\\    \hline
  0.1083   &1.745e-04   &2.019&9.197e-05   &2.087&1.182e-02  &1.009\\    \hline
  0.0541   &4.335e-05   &2.009&2.214e-05   &2.055&5.893e-03  &1.004\\    \hline
  \end{tabular*}
\end{table}

From Table~\ref{Error}, we find that the convergence rates for
$\Bu_h$, $p_h$, $\BB_h$ are given by
 \ben
 \NHonev{\Bu-\Bu_h} \sim O(h^2),\qquad
 \NLtwo{p-p_h} \sim O(h^2),\qquad
 \NHcurl{\BB-\BB_h} \sim O(h).
 \een
Remember that we are using the second-order Lagrangian finite elements for discretizing $\Bu$,
the first-order Lagrangian finite elements for discretizing $p$,
and the first-order N\'{e}d\'{e}lec's edge elements in the second family for discretizing $\BB$. This means that the optimal convergence rates are obtained for all variables.

\begin{example}[Driven Cavity Flow]
The example is the benchmark problem of a driven cavity flow.
The righthand side is given by $\Bf = \mathbf{0}$ and the boundary conditions are given by $\Bg=(g_1, 0, 0)^\top$ and $\BB_s = (1, 0, 0)^\top$ where $g_1=g_1(z)$ is a continuous function and satisfies
 \ben
 g_1 = 1\quad \mathrm{if}~~z = 1; \qquad
 g_1 = 0 \quad \mathrm{if}~~0 \leq z \le h,
 \een
where $h$ is the mesh size. The parameters are set by
 \ben
 R_e = S = 100, \qquad R_m = 1,\qquad \gamma =1.5.
 \een
 \end{example}

The purpose of this example is to verify the effectiveness of the mixed finite element method for engineering benchmark problem
and demonstrate the robustness of the discrete solver with respect
to the mesh size $h$. The computational domain is set by $\Omega = (0,1)^3$. We set the tolerances by $\varepsilon=10^{-6}$ and $\varepsilon_0=10^{-3}$ in Algorithm~\ref{alg:gmres}. The tolerance is $10^{-4}$ for the relative residual of the nonlinear iterations. Table~\ref{cavity:grid} shows the mesh sizes and the numbers of DOFs on the meshes which we are using.
\begin{table}[!htb]
  \centering
  \caption{The mesh sizes and the numbers of DOFs.}
  \label{cavity:grid}
  \begin{tabular*}{10cm}{@{\extracolsep{\fill}}|c|c|c|c|}
  \hline
  Mesh   &$h$    & DOFs for $(\BB,r)$  & DOFs for $(\Bu,p)$    \\\hline
  $\Ct_1$     &0.2165    & 13281     & 15468      \\\hline
  $\Ct_2$     &0.1083    & 97985     & 112724     \\\hline
  $\Ct_3$     &0.0541    & 752001     & 859812        \\\hline
  $\Ct_4$     &0.0271    & 5890817     & 6714692       \\\hline
  \end{tabular*}
\end{table}

\begin{table}[!htb]
  \centering
  \caption{Average GMRES iteration number:$R_e = 100.0, S = 100.0, R_m = 1.0$.}
  \label{cavity:pre}
  \begin{tabular*}{12cm}{@{\extracolsep{\fill}}|c|c|c|}
  \hline
  Mesh        & $N_{\rm picard}\times N_{\rm gmres}$ with $\mathbf{BuBv}$
	  & $N_{\rm picard}\times N_{\rm gmres}$ without $\mathbf{BuBv}$    \\\hline
  $\Ct_1$     & $6\times 51.5 $     & \; \;$7\times 108.1 $     \\\hline
  $\Ct_2$     & $6\times 43.5 $     & \;\; $7\times 102.3 $   \\\hline
  $\Ct_3$     & $6\times 36.8 $     &\;\; $7\times 192.6 $    \\\hline
  $\Ct_4$     & $6\times 32.5 $     &$>7\times 200.0$    \\\hline
  \end{tabular*}
\end{table}

Let $N_{\rm picard}$ denote the number of nonlinear iterations to reduce the relative residual by a factor $10^{-4}$. Let $N_{\rm gmres}$ denote the number of preconditioned GMRES iterations for solving the linearized problem \eqref{Axb}. Therefore, $N_{\rm picard}\times N_{\rm gmres}$ represents the total computational quantity for solving the nonlinear problem \eqref{weakh}.
Remember that the approximate Schur complement $\bbS$ in \eqref{preconditioner} is defined by the bilinear form
 \ben
 \Cf(\Bu_k;\Bu,\Bv) + SR_m (\BB_k\times\Bu,\BB_k\times\Bv) .
 \een
As mentioned in the last paragraph of Subsection~\ref{subsec:coupling}, dropping the second term gives $\bbS=\bbF$. In Table~\ref{cavity:pre}, we show the effectiveness of the preconditioner $\bbP$ with and without the term $\mathbf{BuBv}:=SR_m (\BB_k\times\Bu,\BB_k\times\Bv)$ in $\bbS$. An interesting observation is that, with $\mathbf{BuBv}$, the number of GMRES decays when the mesh is refined successively.
However, without this term, the number of GMRES iterations increases considerably.

\begin{figure}[!htb]
  \centering
  \includegraphics[width=0.3\textwidth]{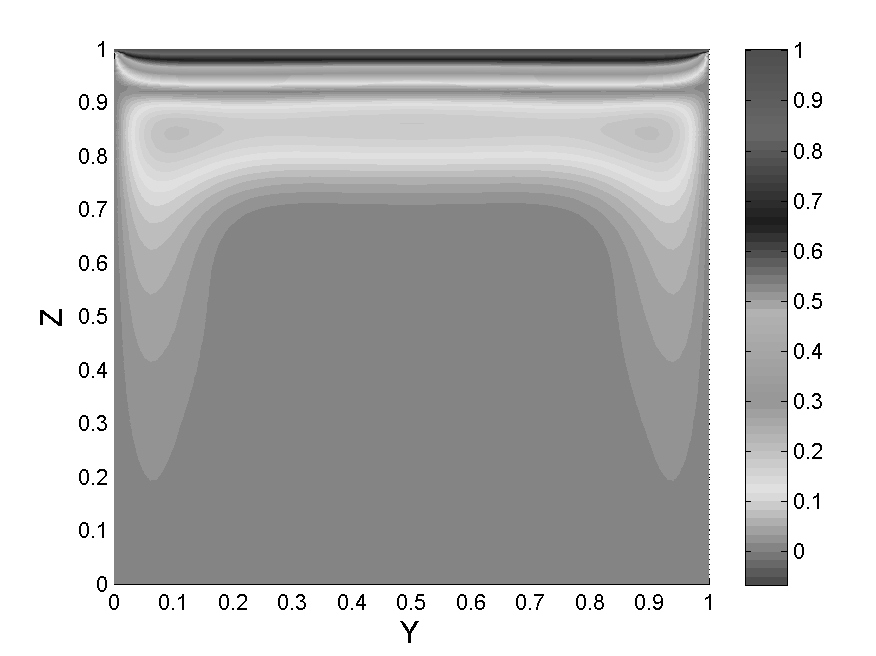}
  \includegraphics[width=0.3\textwidth]{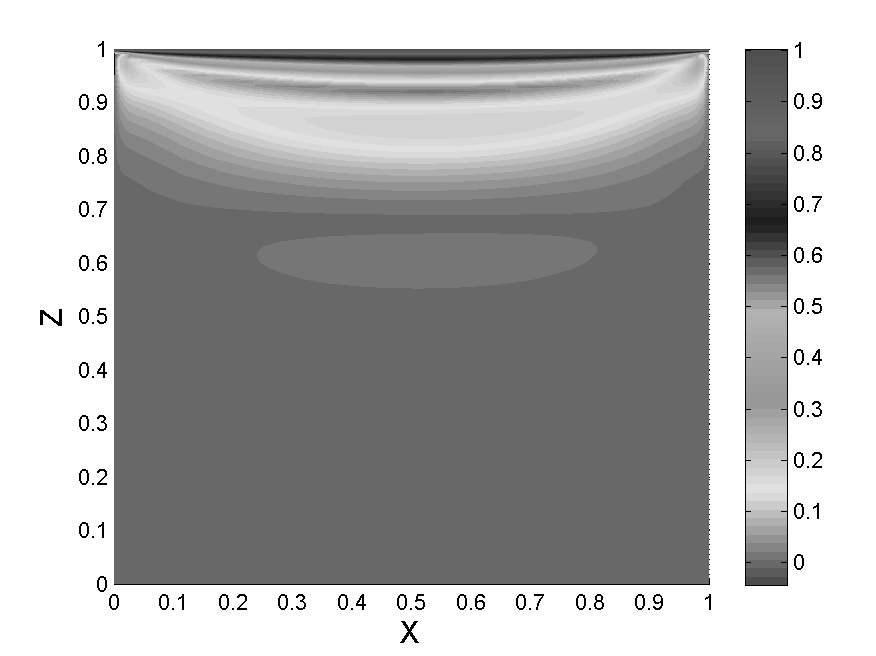}
  \includegraphics[width=0.3\textwidth]{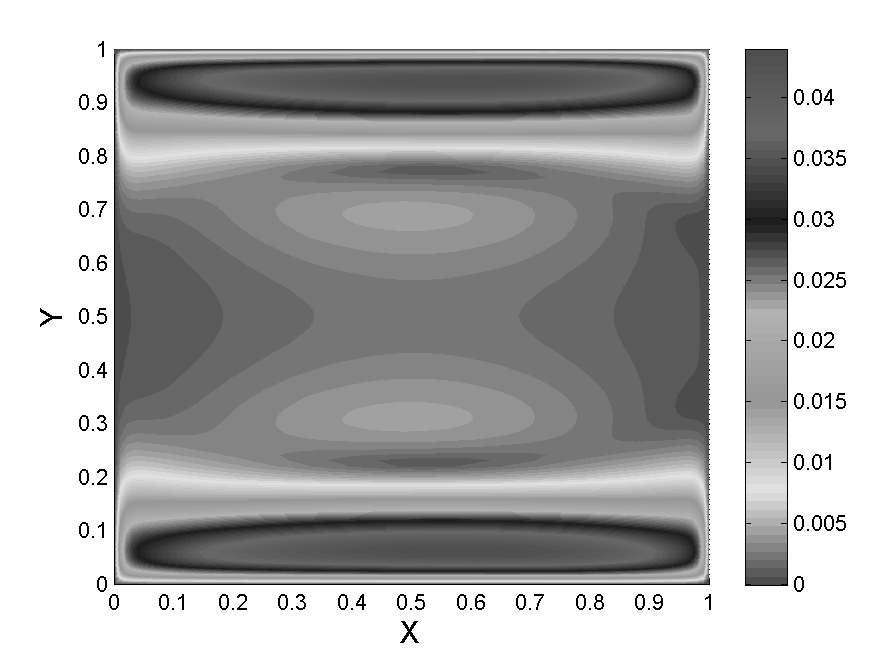}
  \caption{$|\Bu_h|$ on three cross-sections
  $x=0.5$, $y=0.5$, and $z=0.5$ respectively (from left to right).}
  \label{fig:uh}
\end{figure}

\begin{figure}[!htb]
  \centering
  \includegraphics[width=0.3\textwidth]{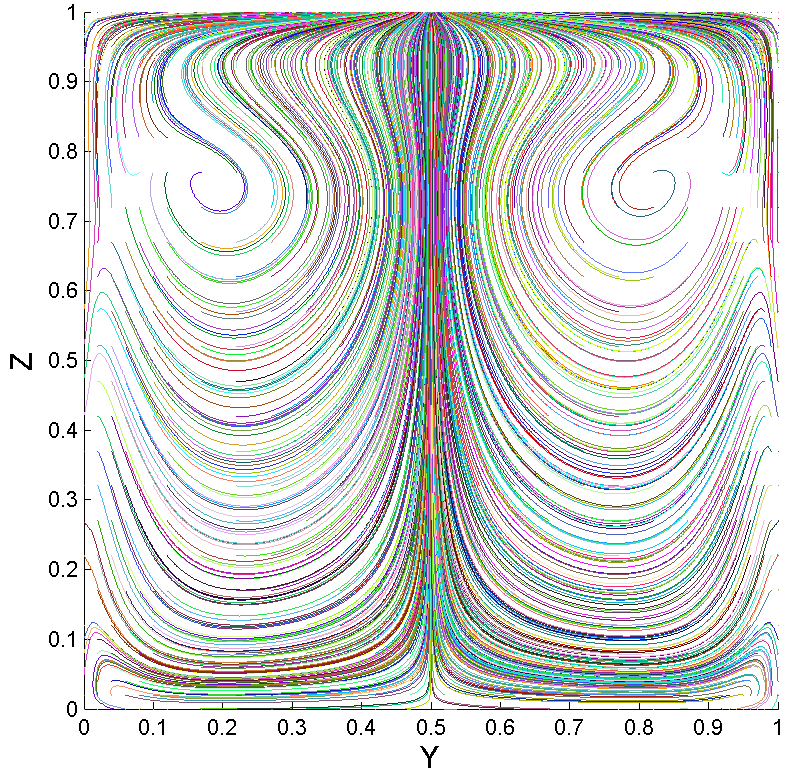}
  \includegraphics[width=0.3\textwidth]{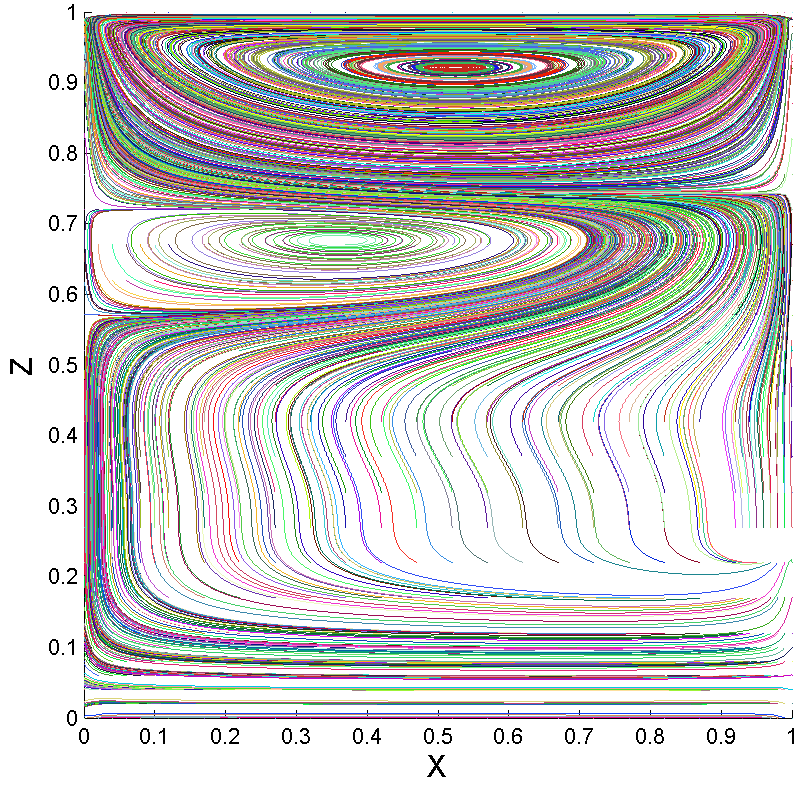}
  \includegraphics[width=0.3\textwidth]{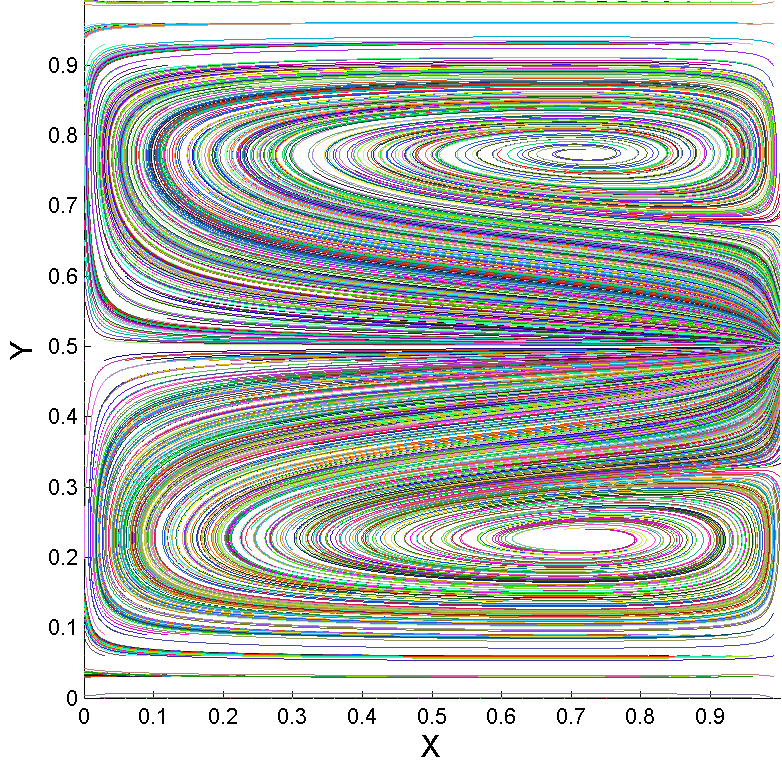}
  \caption{Streamline of the projected velocity.Left:
  $x=0.5$. Middle: $y=0.5$. Right: $z=0.5$.}
  \label{fig:stream}
\end{figure}

\begin{figure}[!htb]
  \centering
  \includegraphics[width=0.3\textwidth]{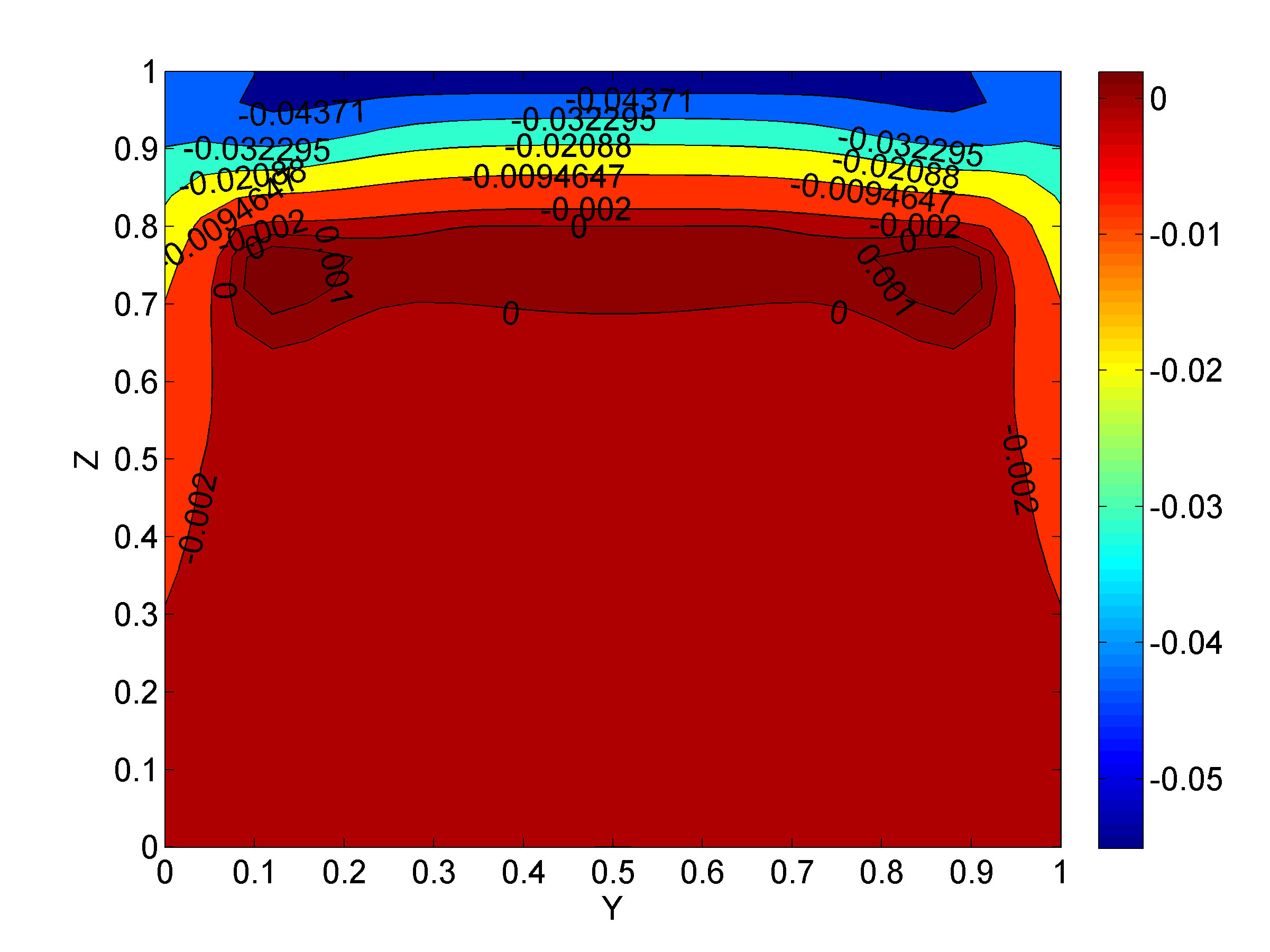}
  \includegraphics[width=0.3\textwidth]{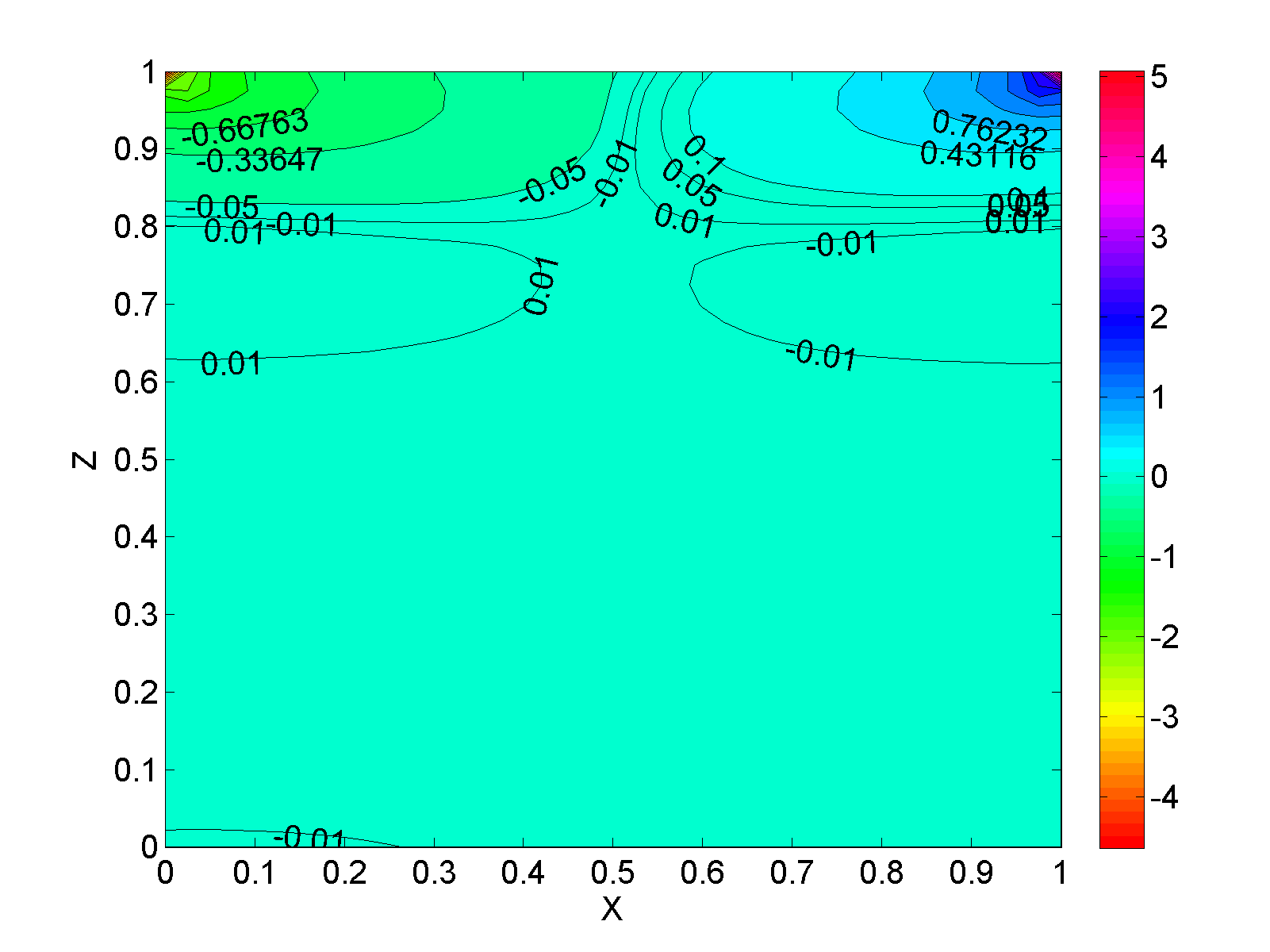}
  \includegraphics[width=0.3\textwidth]{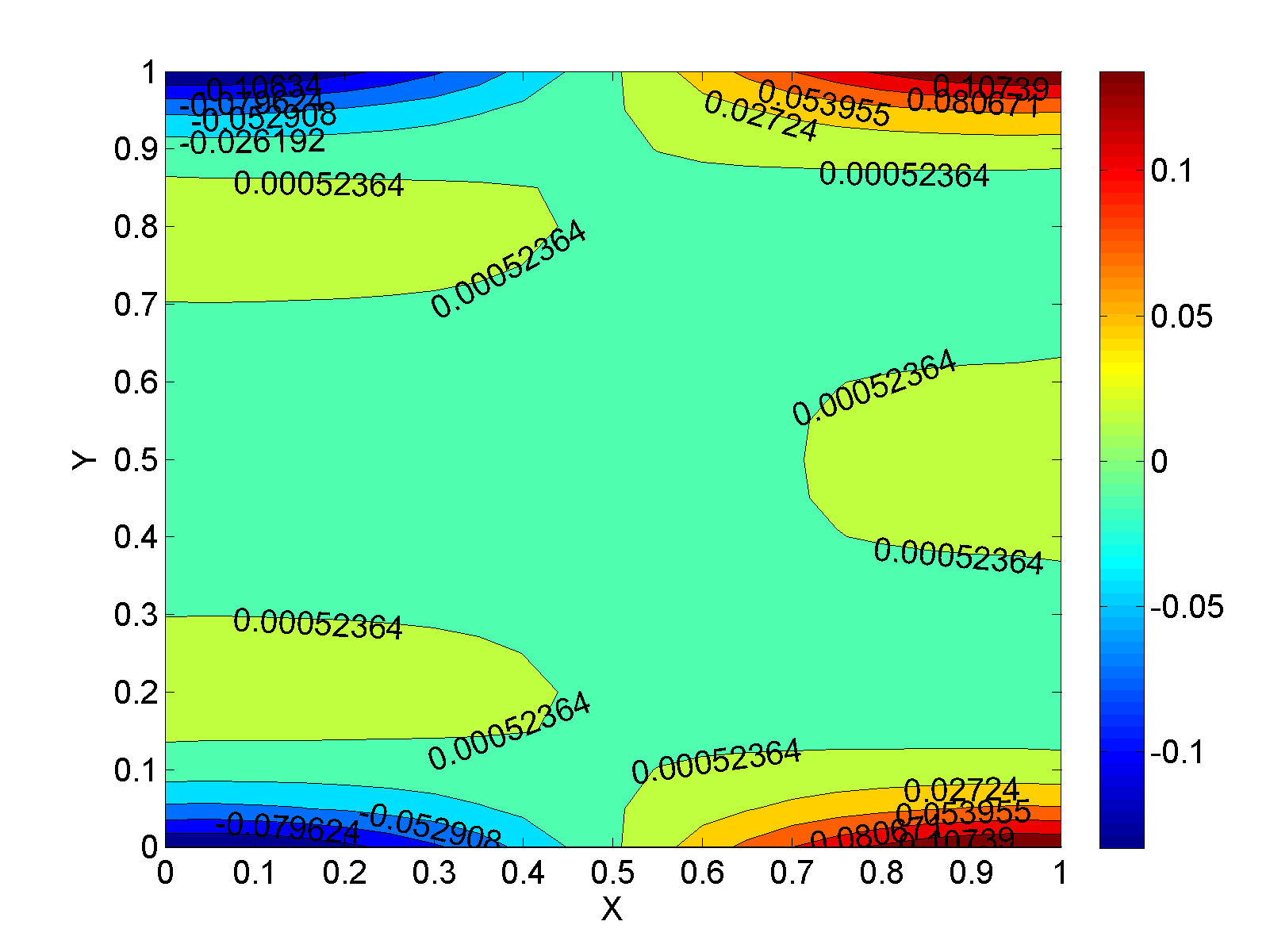}
  \caption{The contour of the pressure $p_h$. Left:
  $x=0.5$. Middle: $y=0.5$. Right: $z=0.5$.}
  \label{fig:ph}
\end{figure}

Now we give the visualization of the simulation results. Firstly, Fig.~\ref{fig:uh} give the grayscale figure of the magnitude of $\Bu_h$. Three 2D projection velocity streamlines of $\Bu_h$ are show
in Fig.~\ref{fig:stream}. Fig.~\ref{fig:ph} shows the contour of the pressure $p_h$.

\begin{figure}[!htb]
  \centering
  \includegraphics[width=\textwidth]{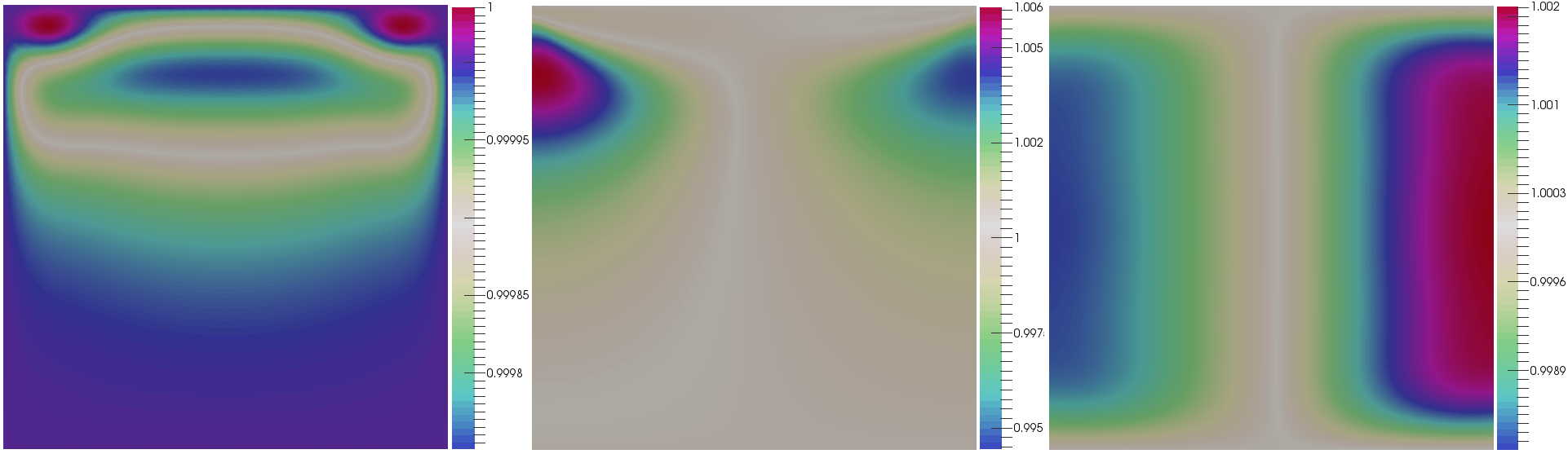}
  \caption{$\SN{\BB_h}$ on three cross-sections
  $x=0.5$, $y=0.5$, and $z=0.5$ respectively (from left to right).}\label{fig:Bh}
\end{figure}

Finally, Fig.~\ref{fig:Bh} shows the distribution of of $\SN{\BB_h}$ on three cross-sections of $\Omega$ at $x=0.5$, $y=0.5$, and $z=0.5$ respectively.

\begin{example}[Scalability]
This example also computes the driven cavity flow and is used to test the scalability of the solver for moderate parameters. The physical parameters $R_e = S = R_m = 1$, and $\gamma = 0.1$.
\end{example}

The error tolerance of the nonlinear iteration is $10^{-4}$ relative to the initial residual. In step 2 of Algorithm~\ref{alg:gmres}, we replace the one level additive Schwarz preconditioner with the BoomerAMG preconditioner \cite{hen02}. This yields better parallel efficiency of the solver for moderate parameters. We carried out the computations on 5 successively refined meshes. Table~\ref{scalability} shows the scalability of the discrete solver by parallel computing on five successively refined meshes. The parallel efficiencies are above $35\%$ and good for such a complex problem.

\begin{table}[!htb]
  \centering
  \caption{Scalability for the discrete solver(Average Iter and Total time).}
  \label{scalability}
  \begin{tabular*}{11cm}{@{\extracolsep{\fill}}|c|c|c|c|c|c|}
  \hline
  Mesh &Total DOFs  & Cores   & $N_{\rm gmres}$ & Time (s)
   & Efficiency \\\hline
  $\Ct_1$   &210709   &1        &16.3  &106.6 & ---\\\hline
  $\Ct_2$   &403221   &2        &14.5  &128.2 &83.2\% \\\hline
  $\Ct_3$   &850965   &4        &15.8  &235.7 &45.2\% \\\hline
  $\Ct_4$   &1611813  &8        &15.8  &300.6 &35.5\% \\\hline
  $\Ct_5$   &3151909  &16       &14.3  &299.7 &35.6\% \\\hline
  \end{tabular*}
\end{table}

For large Reynolds number, the main challenge for the scalability lies in step 2 of Algorithm~\ref{alg:gmres}, that is, the solution of the system of algebraic equations
 \ben
 \bbS\Ve_u = \Vr{(k)} -\bbB^\top\Ve_p.
 \een
We should admit that neither the one level additive Schwarz method nor the classical AMG method can provide ideal scalability solely with our code. Multilevel-based preconditioning and stabilizations for the convection term should be resorted to. This will be our future work and will not be discussed here.


\bibliographystyle{amsplain}

\end{document}